\documentclass[12pt,reqno]{article}
\usepackage{amsmath, amsthm, mathrsfs, graphicx, amsfonts, amssymb, color, cite}
\theoremstyle{plain}

\theoremstyle{Remark}

\newtheorem{prop}{Proposition}[section]
\numberwithin{equation}{section}

\newcommand{\scr}[1]{\mathscr #1}

\newcommand{\set}[1]{\left\{#1\right\}}

\newcommand{\E}{\mathbb{E}}

\renewcommand{\P}{\mathbb P}
\newcommand{\nnb}{\nonumber}

\def\P{\mathbb P}
\def\Z{\mathbb Z}
\def\N{\mathbb N}
\def\B{\scr B}

\def\bg{\begin}
\def\be{\bg{equation}}
\def\de{\end{equation}}

\def\edar{\end{eqnarray}}
\def\beqnn{\begin{eqnarray*}}
\def\eeqnn{\end{eqnarray*}}
\def\lb{\label}
\def\ct{\cite}
\def\l{\left}
\def\r{\right}
\def\fr{\frac}

\def\gm{\gamma}
\def\Gm{\Gamma}
\def\dlt{\delta}

\def\tht{\theta}

\def\kp{\kappa}
\def\lmd{\lambda}

\def\sgm{\sigma}

\def\fa{\forall}

\def\ift{\infty}

\def\rar{\rightarrow}

\def\q{\quad}

\def\var{\text {\rm Var}}

\def\[{\l[} \def\]{\r]}
\def\({\l(} \def\){\r)}
\def\|{\bigg|}
\def\hat{\widehat}
\def\bar{\overline}
\def\tld{\widetilde}

\newcommand{\rf}[1]{(\ref{#1})}

\title{{\bf  Perturbation theory and uniform ergodicity for discrete-time Markov chains}}

\author{
{\bf Yong-Hua Mao$^1$, Yan-Hong Song$^{2, *}$}\\
\footnotesize{$^1$ School of Mathematical Sciences, Beijing Normal University,}\\
\footnotesize{Laboratory of Mathematics and Complex Systems, Ministry of Education,}\\
\footnotesize{Beijing 100875, China}\\
\footnotesize{E-mail:
maoyh@bnu.edu.cn}\\
\footnotesize{$^2$ School of Statistics and Mathematics, Zhongnan University
of Economics and Law,}\\
\footnotesize{Wuhan 430073, China}\\
\footnotesize{E-mail:
$^{*}$songyh@zuel.edu.cn}
}
\date{ }

\begin{document}

\maketitle

\begin{abstract}
We study perturbation theory and uniform ergodicity for discrete-time
Markov chains on general state spaces in terms of the uniform
moments of the first hitting times on some set.
The methods we adopt are different from previous ones.
For reversible and non-negative definite Markov chains,
we first investigate the geometrically ergodic convergence rates.
Based on the estimates, together with a first passage formula,
we then get the convergence rates in uniform ergodicity.
If the transition kernel $P$ is only reversible,
we transfer to study the two-skeleton chain with
the transition kernel $P^2$. At a technical level,
the crucial point is to connect the geometric moments
of the first return times between $P$ and $P^2$.
\end{abstract}

\noindent


%%%%%%%%%%%%%%%%%%%%%%%%%%%%%%%%%%%%%%%%%%%%%%%%%%%%%%%%%%%%%%%%%%%%%%%%%%%%%%%%%%%%
\section{Introduction}

Markov chain Monte Carlo algorithms are
important tools in computational statistics.
The purpose of the algorithm is to draw from
a probability measure $\pi$ by simulating a
Markov chain with transition kernel $P$ such that
$\pi$ is invariant for $P$.
However, it is sometimes impossible to draw from
the transition kernel $P$.
To deal with the difficulty, one may replace $P$
by an approximation $\tld P$.
This leads to a natural question of how small differences
in the transition kernels affect the differences
between their stationary distributions.
Perturbation theory is a common method to study the problem.
It is well-known that there exists an extensive literature
on perturbation bounds for Markov chains.
One group of the results concerns the
sensitivity of uniformly ergodic Markov chains,
see for instance \cite{Kar85a, Kar85b, kar1, Liu12, mit3, RA18}.
The reason is that some practically
important chains are uniformly ergodic,
such as the Metropolis algorithm \ct{Tie94},
some special cases of the Gibbs sampler \ct{RP94, RR98}
and the independent Hastings algorithm \ct{MT96}.

In the paper, we will study perturbation theory and uniform ergodicity
for discrete-time general Markov chains.
Before moving on, let us introduce the basic setup,
the readers are urged to refer \ct{CMF04, MT93}.
Let $\Phi=\{\Phi_n: n\in\Z_+\}$ be a discrete-time homogeneous
Markov chain on a general state space $X$, endowed with
a countably generated $\sigma$-field $\mathscr{B}(X)$.
Denote by $\P_x$ and $\E_x$ the probability and expectation
conditional on $\Phi_0=x$ respectively. Let
$$P^{n}(x, A)=\P_{x}\{\Phi_{n}\in A\},\q
n\in \Z_+, x\in X, A\in\mathscr{B}(X)$$
be the $n$-step transition kernel of the chain,
and it acts on non-negative measurable functions $f$ via
$$P^n f(x)=\int_X f(y)P^n(x, d y),\q n\in\Z_+, x\in X.$$
We assume throughout the paper that the chain $\Phi$
is $\pi$-irreducible for the (unique) invariant probability
measure $\pi$. That is,
if $\pi(A)>0$, there exists $n\in\N$ such that
$P^n(x, A)>0$ for all $x\in X$.
Write $p^n(x, y)$ for the density of $P^n$ with respect to $\pi$,
and $\mathscr{B}^{+}(X)=\l\{A\in\mathscr{B}(X): \pi(A)>0\r\}$
for the sets of positive $\pi$-measure.

For the $\pi$-irreducible chain, it is known (cf. \ct[Chapter 5]{MT93}) that
there always exists some set $A\in\mathscr{B}^+(X)$ satisfying,
for some constants $k, \dlt>0$ and some probability measure $\nu$
on $\mathscr{B}(X)$,
\be\lb{small}
P^k(x, B)\geq\dlt1_A(x)\nu(B),\q x\in X, B\in\mathscr{B}(X).
\de
In what follows, for simplicity of exposition we will consider
the case where $k=1$ in \rf{small}.
That is,
\be\lb{mincondition}
\text{$\exists \dlt>0$ such that $\forall x\in X$,
$B\in\mathscr{B}(X)$,}\q P(x, B)\geq\dlt1_A(x)\nu(B).
\de
\rf{mincondition} is called the minorization condition.
A set $A\in\mathscr{B}(X)$ is called an atom if
\be\lb{atom0}
\text{ $\fa x\in A$, $B\in\mathscr{B}(X)$,}\q P(x, B)=\nu(B).
\de
Moreover, if the atom $A\in\mathscr{B}^{+}(X)$,
then $A$ is called an accessible atom.

For $A\in\mathscr{B}(X)$,
let $\tau_A=\inf\l\{n\geq1: \Phi_n\in A\r\}$ and
$\sgm_A=\inf\l\{n\geq0: \Phi_n\in A\r\}$
be the first return and first hitting times on $A$ respectively.
Denote by
$
F^n(x, A)=\P_x\{\tau_A=n\}
$
the distribution of $\tau_A$, and
%\be\lb{tab}
%{}_AF^n(x, B)=\P_x\set{\Phi_1\notin A,\cdots,
%              \Phi_{n-1}\notin A, \Phi_n\in B}, \ x\in X,B\in\mathscr B(X).
%\de
\be\lb{tab1}
{}_AP^n(x, B)=\P_x\set{\Phi_n\in B, \tau_A\geq n},
\q x\in X,\;A, B\in\mathscr B(X)
\de
the $n$-step taboo probability.
It is clear that $F^n(x, A)={}_AP^n(x, A)$.

By virtue of a discrete-time Phillips' formula
and the coupling technique, we get general results for
perturbation bounds under uniform ergodicity in Section \ref{sec2}.
Phillips' formula enables us to obtain the bounds
for the perturbation of the transition kernels
in uniform total variation norm, and then
the coupling technique helps us to derive the bounds for
the perturbation of the corresponding stationary distributions.
The methods we used here are somewhat
similar to those in \ct{mit1, mit2},
but his methods work well only for Markov chains
on finite state spaces.

The results we obtained in Section \ref{sec2}
are rather extensively applicable once
the uniformly ergodic convergence rates are estimated.
Classically, there are two basic methods to study
convergence rates. One method is to use the renewal theory,
as initiated by Meyn and Tweedie \ct{MT94}.
It requires information about the regeneration
time, which can be obtained by the drift condition.
The other main method, introduced by Rosenthal \ct{Ros95},
is the coupling theory, and relies on estimates of the
coupling time. The methods we adopt in the paper
are different from previous ones.
For reversible and non-negative definite Markov chains
(see Section \ref{ets}),
we first investigate the geometrically ergodic convergence rates.
Our method is in the same spirit of \ct{ST89},
where the geometric convergence rates for Markov chains
on countable state spaces were obtained
via the geometric moments of the first hitting times.
Based on the estimates, together with a first passage formula,
we then get quantitative estimates
on the convergence rates in uniform ergodicity
in terms of the uniform moments of the first hitting times.
When the state space contains an atom,
our results are satisfactory, see Section \ref{atom}.
For non-atomic case,
we use the Nummelin splitting technique to
construct a new Markov chain which admits an atom.
By applying the results in Section \ref{atom} to the
split chain, we obtain the convergence rates
for the original chain, see Section \ref{nonatom}.

Unlike the continuous-time Markov processes,
the discrete-time Markov chains may not
be non-negative definite. For a reversible Markov chain,
we first investigate the two-skeleton chain with
the transition kernel $P^2$, which is also reversible
and always non-negative definite,
and then transfer to $P$, see Section \ref{fr}.
At a technical level, the crucial point is to connect the
geometric moments of the first return times between $P$ and $P^2$.

Finally, in Section \ref{gen}, we study perturbation bound
for general (non-reversible) Markov chains by using a result of
A\"{\i}ssani and Kartashov [1].

%%%%%%%%%%%%%%%%%%%%%%%%%%%%%%%%%%%%%%%%%%%%%%%%%%%%%%%%%%%%%%%%%%%%%%%%%%%%%%%

\section{General results for perturbation bounds}\label{sec2}

In the section, general results for perturbation bounds are obtained.
Recall that the Markov chain $\Phi$ is uniformly ergodic if
$$||P^n-\pi||\rar0,\q n\rar\ift,$$
where
$$||P^n-\pi||=\sup_{x\in X}||P^n(x, \cdot)-\pi||_\var,$$
and for a signed measure $\mu$ on $\mathscr{B}(X)$,
$$||\mu||_\var=\sup_{A\in\mathscr{B}(X)}|\mu|(A)
=\sup_{|f|\leq1}|\mu(f)|,$$
and $\mu(f)=\int_X f d\mu$.
For more details of uniform ergodicity,
see for example \ct{CMF04, CMF05, LZZ08, MYH02, MYH06, MT93}.

Let $\tld P$ be a perturbation of $P$ with
invariant probability measure $\tld\pi$.
For bounded measurable functions $f$,
\be\lb{phillips}
\tld P^n f(x)-P^n f(x)=\sum_{m=0}^{n-1}
\tld P^{n-1-m}(\tld P-P)P^m f(x),\q n\geq 1.
\de
This is a discrete-time version of Phillips' formula in \cite{Phi53}.
In the following, we will study the perturbation bounds
$||\tld P^n-P^n||$ and $||\tld\pi-\pi||_\var$ via the
convergence rates for $P^n$ by the Phillips' formula \rf{phillips}
and the coupling technique.

\bg{thm}\lb{sensitivity}
Let $\set{\gm_n: n\in\Z_+}$ be a positive sequence
such that $\Gm_n:=\sum_{m=0}^{n-1}\gm_m$
with $\Gm_{\ift}:=\sum_{m=0}^{\ift}\gm_m<\ift$.
Assume that
\be\lb{general1}
||P^n-\pi||\leq\gm_n,\q n\geq0.
\de
Then
\be\lb{general2}
||\tld P^n-P^n||\leq\Gm_n||\tld P-P||,
\de
and
\be\lb{general3}
||\tld\pi-\pi||_\var\leq\Gm_{\ift}||\tld P-P||.
\de
\end{thm}
\proof
By the contractivity of $\tld P^n$ and \rf{general1},
for all $0\leq m\leq n-1$ and $|f|\leq1$,
$$\aligned
&\q\sup_{x\in X}\l|\tld P^{n-1-m}(\tld P-P)P^m f(x)\r|
\leq\sup_{x\in X}\l|(\tld P-P)P^m f(x)\r|\\
&=\sup_{x\in X}\l|(\tld P-P)P^m\l(f-\pi(f)\r)(x)\r|
\leq||\tld P-P||\sup_{x\in X}\l|P^m f(x)-\pi(f)\r|\\
&\leq\gm_m||\tld P-P||.
\endaligned
$$
This, together with \rf{phillips}, implies that
for all $|f|\leq1$ and $x\in X$,
\be\lb{stability1}
\l|\tld P^n f(x)-P^n f(x)\r|
\leq\sum_{m=0}^{n-1}\gm_m||\tld P-P||
=\Gm_n||\tld P-P||,
\de
from which we get \rf{general2}.

To pass from \rf{general2} towards \rf{general3},
we use the coupling method.
Since $\pi P^n=\pi$ and $\tld\pi\tld P^n=\tld\pi$, we have
\be\lb{coupling}
\aligned
&\q\l|\tld\pi(f)-\pi(f)\r|
=\l|\int_X \tld P^n f(x)\tld\pi(d x)-\int_X P^n f(y)\pi(d y)\r|\\
&=\l|\int_{X\times X}\l(\tld P^n f(x)-P^n f(y)\r)\tld\pi(d x)\pi(d y)\r|\\
&\leq\int_{X\times X}\l|\tld P^n f(x)-P^n f(y)\r|\tld\pi(d x)\pi(d y)\\
&\leq\int_X\l|\tld P^n f(x)-P^n f(x)\r|\tld\pi(d x)
+\int_{X\times X}\l|P^n f(x)-P^n f(y)\r|\tld\pi(d x)\pi(d y).\nnb
\endaligned
\de
Combining this with \rf{stability1} and noting that
for all $|f|\leq1$ and $x, y\in X$,
$$
\l|P^n f(x)-P^n f(y)\r|
\leq\l|P^n f(x)-\pi(f)\r|+\l|P^n f(y)-\pi(f)\r|\leq2\gm_n,
$$
we get
$$\l|\tld\pi(f)-\pi(f)\r|\leq\Gm_n||\tld P-P||+2\gm_n.$$
Hence
$$\aligned
||\tld\pi-\pi||_\var
&=\sup_{|f|\leq1}\l|\tld\pi(f)-\pi(f)\r|\\
&\leq\lim_{n\rar\ift}\l\{\Gm_n||\tld P-P||+2\gm_n\r\}\\
&=\Gm_{\ift}||\tld P-P||,
\endaligned
$$
which is the desired assertion.
\qed

Note that we can not get \rf{general3} from \rf{general2} directly,
because we only assume $P$ is uniformly ergodic and
the convergence of $\tld P$ is unknown.
Theorem \ref{sensitivity} is rather extensively applicable once
the uniformly ergodic convergence rates are obtained.

\bg{cor}\lb{cor12}
Suppose that the Markov chain $\Phi$ is uniformly ergodic.
That is, there exist constants $\rho<1$ and $C<\ift$ such that
\be\lb{eq211}
||P^n-\pi||\leq C\rho^{n}.\nnb
\de
Then for $n\ge2+\left[\log_{\rho}{(2/C)}\right]$,
\be\lb{ee011}
||\tld P^n-P^n||\leq\left\{2+2\left[\log_{\rho}{(2/C)}\right]+C(1-\rho)^{-1}
\left(\rho^{1+\left[\log_{\rho}{(2/C)}\right]}-\rho^{n}\right)\right\}
||\tld P-P||,\nnb
\de
and
\be\lb{ee02}
||\tld\pi-\pi||_\var\leq \left\{2+2\left[\log_{\rho}{(2/C)}\right]+C(1-\rho)^{-1}
\rho^{1+\left[\log_{\rho}{(2/C)}\right]}\right\}||\tld P-P||.\nnb
\de
\end{cor}

\proof
Since the Markov chain is uniformly ergodic and
$||P^n-\pi||\leq2$,
we can choose $\gm_n$ to be $\gm_n=\min\{2, C\rho^n\}$,
so that for $n\geq2+\left[\log_{\rho}{(2/C)}\right]$,
\be\lb{juti}
\aligned
\Gm_n
&=\sum_{m=0}^{n-1}\gm_m=\sum_{m=0}^{n-1}\min\{2, C\rho^m\}\\
&=\sum_{m=0}^{\left[\log_{\rho}{(2/C)}\right]}2
 +\sum_{m=1+\left[\log_{\rho}^{(2/C)}\right]}^{n-1}C\rho^m\\
&=2+2\left[\log_{\rho}{(2/C)}\right]
 +C(1-\rho)^{-1}\l(\rho^{1+\left[\log_{\rho}{(2/C)}\right]}-\rho^n\r).\nnb\\
\endaligned
\de
Thus,
$$\Gm_{\ift}=2+2\left[\log_{\rho}{(2/C)}\right]+C(1-\rho)^{-1}
\rho^{1+\left[\log_{\rho}{(2/C)}\right]}.$$
Then the desired assertions follow from Theorem \ref{sensitivity}.
\qed

For discrete-time Markov chains,
it is usually convenient to derive the convergence rate by
the Dobrushin's ergodic coefficient.
It is well-known that the chain is uniformly ergodic
if and only if there exists $N\in\N$ such that
\be\lb{ec1}
\dlt:=\fr12\sup_{x,y\in X}||P^N(x,\cdot)-P^N(y,\cdot)||_\var<1.\nnb
\de
Actually, it follows that
$$
||P^n-\pi||\leq2\dlt^{[n/N]},
$$
see e.g. \ct{LW18, mit3, RS18, Sen88}.
The following result gives perturbation bounds
via the Dobrushin's ergodic coefficient, which extends
the previous results for finite and countable Markov chains
in \ct{mit3, Sen88}.

\bg{cor}
Assume that $\dlt<1$. Then
$$||\tld P^n-P^n||\leq\fr{2N(1-\dlt^n)}{1-\dlt}||\tld P-P||,
$$
and
\be\lb{ee121}
||\tld\pi-\pi||_\var\leq\fr{2N}{1-\dlt}||\tld P-P||.\nnb
\de
\end{cor}

%%%%%%%%%%%%%%%%%%%%%%%%%%%%%%%%%%%%%%%%%%%%%%%%%%%%%%%%%%%%%%%%%%%%%%%%%%%%%%%

\section{Uniform ergodicity and perturbation bounds for reversible and
non-negative definite Markov chains}\label{ets}

Recall that the chain is reversible with respect to $\pi$
if
$$\pi(d x)P(x, d y)=\pi(d y)P(y, d x),\q x, y\in X.$$
Since $P^n(x, \cdot)\ll\pi$, we have $p^n(x, y)=p^n(y, x)$
for $\pi\times\pi$-a.s. $(x, y)$ for reversible Markov chains.
Obviously, all Hastings-Metropolis algorithms are
by construction reversible.
If for all $f\in L^2(\pi, X)$,
$$(f, P f)_{L^2(\pi, X)}=\int_{X}f(x)P f(x)\pi(d x)\geq0,$$
the chain is called non-negative definite.
By \ct[Lemma 3.1]{Bax05},
symmetric Metropolis algorithms are non-negative definite.
In the section, we will concentrate on studying
the reversible and non-negative definite Markov chain $\Phi$.

According to \ct[Theorem 16.2.2]{MT93},
the chain $\Phi$ is uniformly ergodic if and only if
\be\lb{hitting}
M:=\sup_{x\in X}\E_x[\sgm_A]<\ift
\de
for some petite set $A\in\mathscr{B}(X)$.
The estimates on the moment of the first hitting time can be obtained by
Foster-Lyapunov drift condition. By \ct[Theorem 11.3.5]{MT93},
we know that if
$$P V\leq V-1+b1_A$$
for some bounded function $V$ and some constant $b<\ift$,
then this gives bounds of the form
$$\E_x[\sgm_A]\leq V(x),\q x\in A^c.$$
In the section, we aim to get quantitative estimates
on the convergence rates in uniform ergodicity
and perturbation bounds in terms of the uniform
moments of the first hitting times $M$.
The method goes as follows.
First, we investigate the geometrically ergodic convergence rates
via the geometric moments of the first return times.
This type of convergence rate has been extensively studied,
see for instance \ct{Bax05, CMF04, CMF05, MS10, MT94, RT99, ST89}
and references therein.
Then, combining the geometric convergence rates with
a first passage formula, we get the convergence rates
in uniform ergodicity.
Based on the results, together with Theorem \ref{sensitivity},
perturbation bounds are finally obtained.
Our study is divided into atomic and non-atomic cases.

\subsection{Atomic case}\label{atom}

It is well-known that the chain $\Phi$ being geometrically ergodic
is equivalent to
\be\lb{geo}
L:=\sup_{x\in A}\E_x[\kp^{\tau_A}]<\ift
\de
for some petite set $A\in\mathscr{B}(X)$ and some constant $\kp>1$,
cf. \ct[Theorem 15.0.1]{MT93}.
By the minimal non-negative solution theory
(see e.g. \ct{CMF04, hg}),
if there exist some function $V\geq1$ and some constant $b<\ift$ satisfying
$$P V\leq\kp^{-1}V+b1_A,$$
then the Foster-Lyapunov drift condition yields a bound of the form
$$\E_x[\kp^{\tau_A}]\leq V(x)+\kp b,\q x\in A.$$
In \ct{ST89}, Sokal and Thomas have studied the geometric
convergence rates by the condition \rf{geo}
for the countable space case.
Since much Markov chain theory on a general state space
can be developed in complete analogy with the countable state
situation if $X$ contains an atom,
we aim to extend part of results in \ct{ST89} to
general chains which admits an atom.

Let ${}_AP$ be a transition kernel from $P$
by restricting on the state space $A^c$.
Denote by $||{}_AP||_{L^2(\pi, A^c)}$ the operator
norm of ${}_AP$ in $L^2(\pi, A^c)$,
and $r_0(P)$ the spectral radium of $P$
in the Hilbert space $\mathbb{H}=\{f\in L^2(\pi, X): \pi(f)=0\}$.
Note that ${}_AP$ is symmetric and non-negative definite. Thus,
\be
(f,{}_AP f)_{L^2(\pi, A^c)}\leq
||{}_AP||_{L^2(\pi, A^c)}||f||_{L^2(\pi, A^c)}^2,
\q f\in {L^2(\pi, A^c)}.\nnb
\de

\bg{lem}\lb{element}
If a reversible and non-negative definite Markov chain
admits an accessible atom, then
$r_0(P)\leq ||{}_AP||_{L^2(\pi, A^c)}$.
\end{lem}

\bg{proof}
Let $f\in\mathbb{H}$
and set $c=\int_A f(x)\pi(d x)/\pi(A)$.
For the accessible atom $A$,
we obtain from the symmetry and \rf{atom0} that
\be\lb{el1}
\aligned
(f, P f)_{L^2(\pi, X)}
&=((f-c1), P(f-c1))_{L^2(\pi, X)}-|c|^2\\
&=\int_{x\in A^c}\int_{y\in A^c}(f(x)-c)(f(y)-c)P(x, d y)\pi(d x)\\
&\q+2\int_{x\in A}\int_{y\in A^c}(f(x)-c)(f(y)-c)P(x, d y)\pi(d x)\\
&\q+\int_{x\in A}\int_{y\in A}(f(x)-c)(f(y)-c)P(x, d y)\pi(d x)-|c|^2\\
&=((f-c1), {}_AP(f-c1))_{L^2(\pi, A^c)}\\
&\q+2\int_{x\in A}\int_{y\in A^c}(f(x)-c)(f(y)-c)\nu(d y)\pi(d x)\\
&\q+\int_{x\in A}\int_{y\in A}(f(x)-c)(f(y)-c)\nu(d y)\pi(d x)-|c|^2\\
&=((f-c1), {}_AP(f-c1))_{L^2(\pi, A^c)}-|c|^2\\
&\leq||{}_AP||_{L^2(\pi, A^c)}||f-c1||_{L^2(\pi, X)}^2-|c|^2\\
&=||{}_AP||_{L^2(\pi, A^c)}(||f||_{L^2(\pi, X)}^2+|c|^2)-|c|^2\\
&\leq||{}_AP||_{L^2(\pi, A^c)}||f||_{L^2(\pi, X)}^2.\nnb\\
\endaligned
\de
Thus, for the non-negative definite Markov chain,
$$r_0(P)=\sup\{(f, P f)_{L^2(\pi, X)}:
\pi(f)=0, ||f||_{L^2(\pi, X)}=1\}\leq||{}_AP||_{L^2(\pi, A^c)},$$
which finishes the proof.
\end{proof}

\bg{prop}\lb{atomerg}
For a reversible and non-negative definite Markov chain,
assume that there exist some accessible atom
$A$ and some constant $\kp>1$ such that \rf{geo} holds.
Then $r_0(P)\leq\kp^{-1}$. Moreover,
$$
\sup_{x\in A}||P^n(x, \cdot)-\pi||_\var\leq\l(\pi(A)^{-1}-1\r)^{1/2}\kp^{-n},
$$
and there exists a constant $C(x)<\ift$ such that
$$
||P^n(x, \cdot)-\pi||_\var\leq C(x)\kp^{-n}, \q\pi\mbox{-a.s.}~x\in X.
$$
\end{prop}

\proof
Note that for all $x\in A^c$,
$$
\l({}_AP^n1\r)(x)=\P_{x}\{\sgm_{A}>n\}
\leq\kp^{-(n+1)}\E_x[\kp^{\sgm_A}].
$$
By this and a generalized form of Kac's formula
(cf. \ct[Lemma 3.4]{Cog75}):
\begin{equation*}
\int_{A}\E_{x}[\kp^{\tau_{A}}]\pi(d x)
=\kp\pi(A)+(\kp-1)\int_{A^c}\E_{x}[\kp^{\sgm_{A}}]\pi(d x),
\end{equation*}
we have for all $f\in L^{\ift}(\pi, A^c)$,
$$
\aligned
&\q\big|(f, {}_AP^n f)_{L^2(\pi, A^c)}\big|
\leq(|f|, {}_AP^n|f|)_{L^2(\pi, A^c)}\\
&\leq ||f||_{\ift}^2\int_{A^c}({}_AP^n1)(x)\pi(d x)\\
&\leq||f||_{\ift}^2\kp^{-(n+1)}
\int_{A^c}\E_{x}[\kp^{\sgm_{A}}]\pi(d x)\\
&\leq||f||_{\ift}^2\kp^{-(n+1)}(\kp-1)^{-1}
\l(\sup_{x\in A}\E_{x}[\kp^{\tau_{A}}]\pi(A)-\kp\pi(A)\r)\\
&\leq C\kp^{-(n+1)}\\
\endaligned
$$
for some $C<\ift$. Since such functions $f$ are dense in $L^2(\pi, A^c)$,
it follows from Lemma \ref{element} and \ct[Proposition 2.5]{ST89}
that
$$
r_0(P)\leq||{}_AP||_{L^2(\pi, A^c)}\leq\kp^{-1}.
$$
Then the spectral mapping theorem yields that
$$||P^n f-\pi(f)||_{L^2(\pi, X)}\leq ||f-\pi(f)||_{L^2(\pi, X)}\kp^{-n}.$$
In particular, from the proof of \ct[Theorem 9.15]{CMF04},
for all probability measure $\mu\ll\pi$,
\be\lb{shang}
||\mu P^n-\pi||_\var\leq\l|\l|\fr{d\mu}{d\pi}-1\r|\r|_{L^2(\pi, X)}\kp^{-n}.
\de

On one hand, set $\mu(dx)=1_A(x)\pi(d x)/\pi(A)$ in \rf{shang}.
Since for the atom $A$,
\be\lb{induction}
P^n(x, \cdot)=\nu(\cdot)\nu^{n-1}(A)+
\sum_{k=1}^{n-1}\int_{A^c}P^k(y, \cdot)\nu(d y)\nu^{n-1-k}(A),\q x\in A\nnb
\de
by induction, which is independent of $ x\in A$,
we get for $x\in A$,
$$
\aligned
||P^n(x, \cdot)-\pi||_\var
&=\sup_{B}\l|\int_X P^n(x, B)1_A(x)\pi(d x)/\pi(A)-\pi(B)\r|\\
&\leq\l(\int\l(\fr{d \mu}{d \pi}\r)^2(x)\pi(d x)-1\r)^{1/2}\kp^{-n}\\
&=\l(\pi(A)^{-1}-1\r)^{1/2}\kp^{-n}.
\endaligned
$$
On the other hand, applying \rf{shang} to $\mu(d y)=P^m(x, d y)$
for $m\leq n$, we get by the reversibility that for $\pi$-a.s. $x\in X$,
\be\lb{vy8}
\aligned
||P^n(x,\cdot)-\pi||_\var&
\leq\l|\l|\fr{d P^m(x, \cdot)}{d\pi}-1\r|\r|_{L^2(\pi, X)}\kp^{-(n-m)}\\
&=\l(\int p^m(x, y)^2\pi(d y)-1\r)^{1/2}\kp^{-(n-m)}\\
&=\l(\int p^m(x, y)p^m(y, x)\pi(d y)-1\r)^{1/2}\kp^{-(n-m)}\\
&=\l[\l(p^{2m}(x, x)-1\r)^{1/2}\kp^{m}\r]\kp^{-n}.\nnb\\
\endaligned
\de
This finishes the proof.
\qed

In order to study the convergence rates in uniform ergodicity,
we still need two lemmas.
The next one connects the uniform geometric moment
of the first return time with the uniform moment of the
first hitting time.

\bg{lem}\lb{hitmoment}
Assume that \rf{hitting} holds for some set $A\in\mathscr{B}(X)$.
Then for all $1<\lmd<e^{1/M}$,
$$
%\sup_{x\in X}\E_x[\lmd^{\sgm_A}]\leq\l(1-M\log\lmd\r)^{-1},\q
\sup_{x\in X}\E_x[\lmd^{\tau_A}]\leq\lmd\l(1-M\log\lmd\r)^{-1}.$$
\end{lem}

\bg{proof}
According to \cite[Theorem 6.3.4]{hg},
$\sup_{x\in X}\E_x[\sgm_A^{\ell}]\leq\ell! M^{\ell}$
for all $\ell\in\N$.
By this and the Taylor expansion of the exponential function,
for all $x\in X$ and $1<\lmd<e^{1/M}$,
\be\lb{m0}
\aligned
\E_x[\lmd^{\sgm_A}]&=\E_x\l[e^{\log\lmd\cdot\sgm_A}\r]
=\sum_{\ell=0}^{\ift}\fr{(\log\lmd)^{\ell}\E_x[\sgm_A^{\ell}]}{\ell!}\\
&\leq\sum_{\ell=0}^{\ift}(\log\lmd)^{\ell}M^{\ell}=\l(1-M\log\lmd\r)^{-1}.\\
\endaligned
\de
Then \ct[Corollary 2.8(1)]{MS14} yields that for $x\in X$,
$$
\aligned
\E_x[\lmd^{\tau_A}]&=\lmd\int_{A^c}\E_y[\lmd^{\sgm_A}]P(x, dy)+\lmd P(x, A)\\
&\leq\lmd\l(1-M\log\lmd\r)^{-1}P(x, A^c)+\lmd P(x, A)\\
&\leq\lmd\l(1-M\log\lmd\r)^{-1}.\\
\endaligned
$$
Thus, we get the desired result.
\end{proof}

By using the following first passage formula,
Proposition \ref{atomerg} can be applied to
get the convergence rate in uniform ergodicity.

\bg{lem}\lb{renewal}
Let $A\in\mathscr{B}(X)$. For all $x\in X$,
\be\lb{passage}
||P^n(x, \cdot)-\pi||_\var\leq 2\P_x\{\tau_A\geq n+1\}
+\sum_{m=1}^n\sup_{y\in A}||P^{n-m}(y, \cdot)-\pi||_\var~F^m(x, A).
\de
\end{lem}

\proof
For all $x\in X$ and $B\in\mathscr{B}(X)$,
the following decomposition formula holds
by using the taboo probability (cf. \rf{tab1}):
$$
\aligned
P^n(x, B)&=\P_x\l\{\Phi_n\in B, \tau_A\geq n+1\r\}
+\P_x\l\{\Phi_n\in B, \tau_A\leq n\r\}\\
&=\P_x\l\{\Phi_n\in B, \tau_A\geq n+1\r\}
+\sum_{m=1}^n\int_A P^{n-m}(y, B){}_AP^m(x, d y).
\endaligned
$$
It follows that
$$
\aligned
P^n(x, B)-\pi(B)
=&\P_x\l\{\Phi_n\in B, \tau_A\geq n+1\r\}-\pi(B)\P_x\l\{\tau_A\geq n+1\r\}\\
&+\sum_{m=1}^n\int_A \l(P^{n-m}(y, B)-\pi(B)\r){}_AP^m(x, dy).\\
\endaligned
$$
Therefore,
$$
\aligned
||P^n(x, \cdot)-\pi||_\var
&\leq\sup_{B\in \B(X)}\Bigg\{\P_x\l\{\Phi_n\in B, \tau_A\geq n+1\r\}
 +\pi(B)\P_x\l\{\tau_A\geq n+1\r\}\\
&\q\q\q\q\q+\sum_{m=1}^n\int_A \l|P^{n-m}(y, B)-\pi(B)\r|{}_AP^m(x, d y)\Bigg\}\\
&\le2\P_x\{\tau_A\geq n+1\}
+\sum_{m=1}^n\sup_{y\in A}||P^{n-m}(y, \cdot)-\pi||_\var~F^m(x, A),\\
\endaligned
$$
where the last inequality holds since ${}_AP^m(x, A)=F^m(x, A)$.
\qed

\bg{thm}\lb{key1}
For a reversible and non-negative definite Markov chain,
assume that there exists some accessible atom
$A$ such that \rf{hitting} holds.
Then for all $1<\lmd<e^{1/M}$,
$$\l|\l|P^n-\pi\r|\r|\leq D_1e^{-n/M}+E_1\lmd^{-n},$$
where
\be\lb{changshu1}
\aligned
&D_1=C_1\l(1-\fr{e^{1/M}-1}{e^{1/M}-\lmd}M_1\r),\q
E_1=M_1\l((2-C_1)^+\lmd^{-1}+\fr{e^{1/M}-1}{e^{1/M}-\lmd}C_1\r),\\
&C_1=\l(\pi(A)^{-1}-1\r)^{1/2},\q M_1=\lmd(1-M\log\lmd)^{-1}.
\endaligned
\de
\end{thm}

\proof
By Lemma \ref{hitmoment}, $\sup_{x\in X}\E_x[\lmd^{\tau_A}]<\ift$
for all $1<\lmd<e^{1/M}$.
Then Proposition \ref{atomerg} implies that $r_0(P)\leq\lmd^{-1}$
for all $1<\lmd<e^{1/M}$, so that $r_0(P)\leq e^{-1/M}$. Thus,
\be\lb{ui}
\sup_{x\in A}||P^n(x, \cdot)-\pi||_\var\leq\l(\pi(A)^{-1}-1\r)^{1/2}e^{-n/M}.
\de
According to Lemma \ref{renewal} and \rf{ui},
we have for all $x\in X$,
\be\lb{3361}
\aligned
||P^n(x, \cdot)-\pi||_\var
&\leq 2\P_x\{\tau_A\geq n+1\}
+C_1\sum_{m=1}^n e^{-(n-m)/M}\P_x\l\{\tau_A=m\r\}.\\
\endaligned
\de
Let $a_n=\P_x\{\tau_A\geq n\}$ and
$b_n=e^{-n/M}$ for $n\geq0$.
It follows from Abel's theorem that
\be\lb{3362}
\aligned
&\q\sum_{m=1}^n e^{-(n-m)/M}\P_x\l\{\tau_A=m\r\}
=\sum_{m=1}^n b_{n-m}(a_m-a_{m+1})\\
&=b_n-a_{n+1}+\sum_{m=1}^n(b_{n-m}-b_{n+1-m})a_m\\
&=e^{-n/M}-\P_x\{\tau_A\geq n+1\}+(1-e^{-1/M})e^{-n/M}
\sum_{m=1}^n e^{m/M}\P_x\{\tau_A\geq m\}.\\
\endaligned
\de
Combining \rf{3361} with \rf{3362},
we obtain from Lemma \ref{hitmoment} that for all $1<\lmd<e^{1/M}$,
\be\lb{363}
\aligned
||P^n-\pi||
&\leq C_1e^{-n/M}+(2-C_1)^+\sup_{x\in X}\P_x\{\tau_A\geq n+1\}\\
&\q+C_1(1-e^{-1/M})e^{-n/M}
\sum_{m=1}^n e^{m/M}\sup_{x\in X}\P_x\{\tau_A\geq m\}\\
&\leq  C_1e^{-n/M}+(2-C_1)^+\lmd^{-(n+1)}\sup_{x\in X}\E_x\l[\lmd^{\tau_A}\r]\\
&\q+C_1(1-e^{-1/M})e^{-n/M}\sum_{m=1}^n
 \l(e^{1/M}/\lmd\r)^m\sup_{x\in X}\E_x\l[\lmd^{\tau_A}\r]\\
&\leq C_1\l(1-\fr{e^{1/M}-1}{e^{1/M}-\lmd}M_1\r)e^{-n/M}
 +M_1\l((2-C_1)^+\lmd^{-1}+\fr{e^{1/M}-1}{e^{1/M}-\lmd}C_1\r)\lmd^{-n},\nnb
\endaligned
\de
which finishes the proof.
\qed

Combining Theorems \ref{sensitivity} with \ref{key1},
we obtain the following perturbation results directly.

\bg{thm}
Under assumptions of Theorem \ref{key1}, we have for all $1<\lmd<e^{1/M}$,
\be\lb{general2a}
\aligned
||\tld P^n-P^n||
\leq\l\{\fr{e^{1/M}}{e^{1/M}-1}D_1(1-e^{-n/M})+
\fr{\lmd}{\lmd-1}E_1(1-\lmd^{-n})\r\}||\tld P-P||,\nnb\\
\endaligned
\de
and
\be\lb{general3a}
\aligned
||\tld\pi-\pi||_\var\leq\l\{\fr{e^{1/M}}{e^{1/M}-1}D_1+\fr{\lmd}{\lmd-1}E_1\r\}
||\tld P-P||,\nnb\\
\endaligned
\de
where $D_1$ and $E_1$ are defined in \rf{changshu1}.
\end{thm}

%%%%%%%%%%%%%%%%%%%%%%%%%%%%%%%%%%%%%%%%%%%%%%%%%%%%%%%%%%%%%%%%%%%%%%%%%%%%%%%
\subsection{Non-atomic case}\label{nonatom}

On general state spaces, however, accessible atoms are less frequent.
Fortunately, by suitably extending the
probabilistic structure of the $\pi$-irreducible chain, we can artificially
construct a new Markov chain which contains an atom, and this
allows much of the critical analysis to follow the form
of atom case. For more details,
one can refer to Athreya and Ney \ct{AN80} and Nummelin \ct{Num78, Num84}.
In the section, we will focus on the Nummelin splitting.

Suppose that the minorization condition \rf{mincondition}
holds for some set $A\in\mathscr{B}^+(X)$.
Then the construction can be carried out by splitting
the state space $X$, the measure on
$\mathscr{B}(X)$ and the transition kernel $P$ separately.
First, we split the state space $X$ by writing
$\check{X}=X\times\{0, 1\}$, where $X_0=X\times\{0\}$
and $X_1=X\times\{1\}$ are equipped with $\sgm$-field
$\mathscr{B}(X_0)$ and $\mathscr{B}(X_1)$ respectively.
Let $\mathscr{B}(\check{X})$ be the $\sgm$-field of the
subsets of $\check{X}$ generated by $\mathscr{B}(X_0)$
and $\mathscr{B}(X_1)$. Write $x_0\in X_0$ and $x_1\in X_1$ for the two ``copies of $x$" and
$B_0\subseteq X_0$ and $B_1\subseteq X_1$
for the ``copies of $B$" for $B\in\mathscr{B}(X)$.
Let $\mu$ be any measure on $\mathscr{B}(X)$. We next split the
measure $\mu$ into two measures by defining the measure $\mu^*$
on $\mathscr{B}(\check{X})$ through
\be\lb{measure}
\aligned
\begin{array}{ll}
\mu^*(B_0)=(1-\dlt)\mu(A\cap B)+\mu(A^c\cap B);& \\
\mu^*(B_1)=\dlt\mu(A\cap B),& \\
\end{array}\\
\endaligned
\de
where $\dlt$ and $A$ are the constant and the set respectively
in \rf{mincondition}.
The important point to notice is $\mu$ is the marginal measure of $\mu^*$
in the sense that
\be\lb{bianmiu}
\mu^*(B_0\cup B_1)=\mu(B),\q B\in\mathscr{B}(X).
\de
Finally, define the split transition kernel $\check{P}(x_i, \cdot)$
for $x_i\in\check{X}$ by
\be\lb{kernel}
\aligned
\begin{array}{ll}
\check{P}(x_0, \cdot)=P(x, \cdot)^*,& x_0\in X_0\setminus A_0;\\
\check{P}(x_0, \cdot)=(1-\dlt)^{-1}\l(P(x, \cdot)^*-\dlt\nu^*(\cdot)\r),
& x_0\in A_0;\\
\check{P}(x_1, \cdot)=\nu^*(\cdot),& x_1\in X_1,\nnb\\
\end{array}\\
\endaligned
\de
where $\dlt$, $A$ and $\nu$ are the constant, the set
and the measure respectively in \rf{mincondition}.
The minorization condition ensures that the transition kernel
$\check{P}$ is meaningful.
According to the above three steps, we get a split chain
$\check{\Phi}$ on $\l(\check{X}, \mathscr{B}(\check{X})\r)$
with transition kernel $\check{P}$,
which admits an atom $X_1$.
Since $\check{P}^n(x_i, X_1\setminus A_1)=0$
for all $n\geq1$ and $x_i\in\check{X}$, $A_1$
is the set which is reached with positive probability. Hence,
we denote by $A_1$ the atom of the split chain $\check{\Phi}$.
It is clear from \rf{measure} that $\pi^*(A_1)=\dlt\pi(A)$.

The splitting technique is important because of
the various properties that $\check{\Phi}$
inherits from, or passes on to, $\Phi$.
From \ct[Theorem 5.1.3 and Proposition 10.4.1]{MT93},
we get the following proposition.

\begin{prop}\lb{connect}
$(i)$
The chain $\Phi$ is the marginal chain of $\check{\Phi}$.
That is, for all initial distribution $\mu$ on $\mathscr{B}(X)$
and all $B\in\mathscr{B}(X)$,
$$\int_X P^n(x, B)\mu(d x)=
\int_{\check{X}}\check{P}^n(x_i, B_0\cup B_1)\mu^*(d x_i).$$

$(ii)$
If the chain $\Phi$ is $\psi$-irreducible with $\psi(A)>0$,
then $\check{\Phi}$ is $\nu^*$-irreducible;
and the chain $\Phi$ is $\psi$-irreducible
if $\check{\Phi}$ is $\psi^*$-irreducible.

$(iii)$
If the measure $\pi$ is invariant for $\Phi$,
then $\pi^*$ is invariant for $\check{\Phi}$;
and if the measure $\check{\pi}$ is invariant for $\check{\Phi}$,
then the measure $\pi$ on $\mathscr{B}(X)$ defined by
$$\pi(B)=\check{\pi}(B_0\cup B_1),\q B\in\mathscr{B}(X)$$
is invariant for $\Phi$, and $\check{\pi}=\pi^*$.
\end{prop}

Let $\check{\P}_{x_i}$ and $\check{\E}_{x_i}$
be the probability and expectation for the split
chain started with $\check{\Phi}_0=x_i$ respectively.
Since $A_1$ is the atom of $\check{\Phi}$,
we write for simplicity $\check{\P}_{A_1}=\check{\P}_{x_1}$
and $\check{\E}_{A_1}=\check{\E}_{x_1}$ for $x\in A$.
Define
\be\lb{expectation}
\check{\E}_x=\l(1-\dlt 1_A(x)\r)\check{\E}_{x_0}
+\dlt 1_A(x)\check{\E}_{x_1}.
\de
Clearly, $\check{\E}_x$ agree with $\E_x$ on $\mathscr{B}(X)$.
Let
$\check{\tau}_{A_1}=\inf\{n\geq1: \check{\Phi}_n\in A_1\}$
be the first return time to $A_1$ for $\check{\Phi}$. Denote by
$\check{F}^n(x_i, A_1)=\check{\P}_{x_i}\{\check{\tau}_{A_1}=n\}$
the distribution of $\check{\tau}_{A_1}$.
By \ct[LEMMA A.1]{Bax05} and \ct[Corollary 2.8]{MS14},
we obtain the following lemma, which shows the relationship
of the geometric moments between $\check\tau_{A_1}$ and $\tau_A$.

\bg{lem}\lb{bax}
Assume that \rf{mincondition} and \rf{geo} hold
for some set $A\in\mathscr{B}(X)$ and some constant $\kp>1$. Then
\be\lb{ccc}
\check{\E}_{x_i}[\lmd^{\check{\tau}_{A_1}}]\leq
\fr{\dlt\check{\E}_{x_i}[\lmd^{\tau_{A}}]}
{1-(1-\dlt)\sup\limits_{x\in A}\check{\E}_{x_0}[\lmd^{\tau_{A}}]},\nnb
\q x_i\in \check{X}
\de
for all $1<\lmd<\kp$ such that
$(1-\dlt)\sup_{x\in A}\check{\E}_{x_0}[\lmd^{\tau_{A}}]<1$.
\end{lem}

According to the above lemma,
we estimate the geometric moment of $\check\tau_{A_1}$.

\bg{lem}\lb{baxmore}
$(i)$
Assume that \rf{mincondition} and \rf{geo} hold
for some set $A\in\mathscr{B}(X)$ and some constant $\kp>1$. Then
for all $1<\lmd<\kp\wedge(1-\dlt)^{-1/\alpha}$,
\be\lb{mki}
\check{\E}_{A_1}[\lmd^{\check{\tau}_{A_1}}]
\leq\fr{\dlt\lmd^{\beta}}{1-(1-\dlt)\lmd^{\alpha}},\nnb
\de
where
\be\lb{alpha}
\alpha=\l(\log\fr{L-\dlt\kp}{1-\dlt}\r)\big/\l(\log\kp\r),\q
\beta=\l(\log\fr{L-(1-\dlt)\kp}{\dlt}\r)\big/\l(\log\kp\r).
\de

$(ii)$
Assume that \rf{mincondition} and \rf{hitting} hold
for some set $A\in\mathscr{B}(X)$.
Then for $1<\lmd<e^{1/M}$ satisfying
\be\lb{lmd1}
\lambda<(1+\delta\lambda)({1-M\log\lmd}),
\de
we have
\be\lb{buchong}
\sup_{x_i\in\check{X}}\check{\E}_{x_i}[\lmd^{\check{\tau}_{A_1}}]
\leq\frac{ {\delta\lambda}}{(1+\delta \lambda)({1-M \log\lmd })-\lambda}.\nnb
\de

\end{lem}

\proof
(i)
From Lemma \ref{bax}, we need to estimate
$\sup_{x\in A}\check{\E}_{x_0}[\lmd^{\tau_{A}}]$
and $\check{\E}_{A_1}[\lmd^{\tau_{A}}]$ separately.
By Jensen's inequality and \rf{expectation},
we get for all $1<\lmd<\kp$,
\be\lb{expection2d}
\aligned
&\q\sup_{x\in A}\check{\E}_{x_0}[\lmd^{\tau_{A}}]
\leq\sup_{x\in A}
 \l(\check{\E}_{x_0}[\kp^{\tau_{A}}]\r)^{(\log\lmd)/(\log\kp)}\\
&=\sup_{x\in A}\l(\fr{\E_x[\kp^{\tau_A}]-\dlt\check{\E}_{A_1}
 [\kp^{\tau_A}]}{1-\dlt}\r)^{(\log\lmd)/(\log\kp)}\\
&\leq\l(\fr{L-\dlt\kp}{1-\dlt}\r)^{(\log\lmd)/(\log\kp)}\\
&=\lmd^{\alpha}.\\
\endaligned
\de
Noting that
\be\lb{nui}
\aligned
\check{\E}_{A_1}[\kp^{\tau_{A}}]
&\leq\dlt^{-1}\sup_{x\in A}\l(\E_x[\kp^{\tau_A}]
 -(1-\dlt)\check{\E}_{x_0}[\kp^{\tau_A}]\r)
\leq\fr{L-(1-\dlt)\kp}{\dlt},\nnb\\
\endaligned
\de
so similarly, for all $1<\lmd<\kp$,
\be\lb{expections}
\aligned
\check{\E}_{A_1}[\lmd^{\tau_{A}}]
\leq\l(\fr{L-(1-\dlt)\kp}{\dlt}\r)^{(\log\lmd)/(\log\kp)}=\lmd^{\beta}.
\endaligned
\de
Thus, the desired assertion holds by Lemma \ref{bax}
with \rf{expection2d} and \rf{expections}
for all $1<\lmd<\kp\wedge(1-\dlt)^{-1/\alpha}$.

(ii)
Combining \rf{expectation} with Lemma \ref{hitmoment},
for all $1<\lmd<e^{1/M}$ satisfying \rf{lmd1},
\be\lb{expection2}
\aligned
(1-\dlt)\sup_{x\in A}\check{\E}_{x_0}[\lmd&^{\tau_{A}}]
=\sup_{x\in A}{\E_x[\lmd^{\tau_A}]-\dlt\check{\E}_{A_1}
 [\lmd^{\tau_A}]}\\
&\leq{\lmd(1-M\log\lmd)^{-1}-\dlt\lmd}<1.\nnb\\
\endaligned
\de
Thus, by Lemma \ref{bax} and noting that
$$
\sup_{x_i\in\check{X}}\check{\E}_{x_i}[\lmd^{\tau_{A}}]=
\sup_{x\in X}\E_{x}[\lmd^{\tau_{A}}]\leq \lmd(1-M\log \lmd)^{-1},
$$
we get the desired result.
\qed

\bg{rem}
It should be pointed out that \rf{lmd1} is meaningful, since
$$
f(\lmd):=(1+\delta \lambda)({1-M \log\lmd })/\lmd
$$
is decreasing with $f(1)=1+\dlt>1$ and $f(e^{1/M})=0$.
Moreover, note that $\log\lmd< \lmd-1$ for all $\lmd>1$,
we can solve the following quadric inequality
$$
\dlt M\lmd^2+(M+1)(1-\dlt)\lmd-(M+1)\leq0
$$
to get an estimation of $\lmd$.
\end{rem}

Applying the techniques used in Proposition \ref{atomerg}
to the split chain $\check{\Phi}$, we get the following result,
which is part of counterparts to Proposition \ref{atomerg},
and is also important out scope of the paper.

\bg{prop}\lb{minerg-prop}
For a reversible and non-negative definite Markov chain,
assume that there exist some set $A\in\mathscr{B}^+(X)$
and some constant $\kp>1$ such that
\rf{mincondition} and \rf{geo} hold.
Then there exists a constant $C(x)<\ift$ such that
\be\lb{dak}
||P^n(x, \cdot)-\pi||_\var\leq C(x)K^{-n}, \q\pi\mbox{-a.s.}~x\in X,\nnb
\de
where $K=\kp\wedge(1-\dlt)^{-1/\alpha}$
and $\alpha$ is defined in \rf{alpha}.
\end{prop}

\bg{rem}
Noting that $\lim_{\dlt\rar1}(1-\dlt)^{-1/\alpha}=\kp$,
so $K=\kp$ for atomic case.
That is, the convergence rate in Proposition \ref{minerg-prop}
is consistent with that in Proposition \ref{atomerg} when
the set $A$ is an accessible atom.
\end{rem}

\proof
According to Proposition \ref{connect},
the split chain $\check{\Phi}$ possessing an accessible atom $A_1$
is reversible with respect to $\pi^*$
and non-negative definite.
Thus, from Lemma \ref{baxmore}(i) and a similar proof as that of
Proposition \ref{atomerg}, we have for $\mu^*\ll\pi^*$,
\be\lb{split-con}
||\mu^* \check{P}^n-\pi^*||_\var
\leq\l|\l|\fr{d\mu^*}{d\pi^*}-1\r|\r|_{L^2(\pi^*, \check{X})}K^{-n}.
\de
In the following, we transfer the above result to the original chain $\Phi$.
By Proposition \ref{connect}(i) and \rf{bianmiu},
\be\lb{qian}
\aligned
||\mu P^n-\pi||_\var
&=\sup_{B\in\mathscr{B}(X)}\bigg|\int_X P^n(x, B)\mu(d x)-\pi(B)\bigg|\\
&=\sup_{B\in\mathscr{B}(X)}\bigg|\int_{\check X}
  \check{P}^n(x_i, B_0\cup B_1)\mu^*(d x_i)-\pi^*(B_0\cup B_1)\bigg|\\
&\leq\sup_{\check{B}\in\mathscr{B}(\check{X})}\bigg|\int_{\check X}
 \check{P}^n(x_i, \check{B})\mu^*(d x_i)-\pi^*(\check{B})\bigg|\\
&=||\mu^* \check{P}^n-\pi^*||_\var,\\
\endaligned
\de
and by \rf{measure},
\be\lb{zhong}
\aligned
&\q\l|\l|\fr{d\mu^*}{d\pi^*}-1\r|\r|_{L^2(\pi^*, \check{X})}^2
=\int_{X_0\cup X_1}\fr{\mu^*(d y_i)^2}{\pi^*(d y_i)}-1\\
&=\int_{X}\fr{((1-\dlt)\mu(A\cap d y)+\mu(A^c\cap d y))^2}
  {(1-\dlt)\pi(A\cap d y)+\pi(A^c\cap d y)}
  +\int_{X}\fr{(\dlt\mu(A\cap d y))^2}{\dlt\pi(A\cap d y)}-1.\\
\endaligned
\de
Combining \rf{split-con} with \rf{qian} and \rf{zhong}, we get
\be\lb{hou}
\aligned
||\mu P^n-\pi||_\var
&\leq\Big(\int_{X}\fr{((1-\dlt)\mu(A\cap d y)+\mu(A^c\cap d y))^2}
  {(1-\dlt)\pi(A\cap d y)+\pi(A^c\cap d y)}\\
&\q\q+\int_{X}\fr{(\dlt\mu(A\cap d y))^2}{\dlt\pi(A\cap d y)}-1\Big)^{1/2}K^{-n}.\nnb\\
\endaligned
\de
Applying the above inequality to $\mu(d y)=P^m(x, d y)$
for $m\leq n$, we have by the reversibility that for $\pi$-a.s. $x\in X$,
\be\lb{liang1}
\aligned
||P^n(x,\cdot)-\pi||_\var
&\leq\Big(\int_{X}\fr{((1-\dlt)P^m(x, A\cap d y)+P^m(x, A^c\cap d y))^2}
  {(1-\dlt)\pi(A\cap d y)+\pi(A^c\cap d y)}\\
&\q\q+\int_{X}\fr{(\dlt P^m(x, A\cap d y))^2}{\dlt\pi(A\cap d y)}-1\Big)^{1/2}K^{-(n-m)}\\
&=\Big((1-\dlt)\int_A p^m(x, y)^2\pi(d y)+
\int_{A^c} p^m(x, y)^2\pi(d y)\\
&\q\q+\dlt\int_A p^m(x, y)^2\pi(d y)-1\Big)^{1/2}K^{-(n-m)}\\
&=\l(\int p^m(x, y)^2\pi(d y)-1\r)^{1/2}K^{-(n-m)}\\
&=\l[\l(p^{2m}(x, x)-1\r)^{1/2}K^m\r]K^{-n}.\nnb\\
\endaligned
\de
Thus, the proof is finished.
\qed

We can now move from the geometric ergodicity result to
uniform ergodicity and perturbation bounds for non-atomic case.

\bg{thm}\lb{minerg}
For a reversible and non-negative definite Markov chain,
assume that there exists a set $A\in\mathscr{B}^+(X)$ such that
\rf{mincondition} and \rf{hitting} hold.
Then for all $1<\lmd<e^{1/M}$ satisfying \rf{lmd1},
\be\lb{chang}
\aligned
\l|\l|P^n-\pi\r|\r|\leq\l(D_2+E_2 n\r)\lmd^{-n},\nnb\\
\endaligned
\de
where
\be\lb{define-b}
\aligned
&D_2=C_2+(2-C_2)^+\lmd^{-1}M_2,\q E_2=C_2(1-\lmd^{-1})M_2,\\
&C_2=\l((\dlt\pi(A))^{-1}-1\r)^{1/2},\q
M_2=\frac{{\delta\lmd}}{(1+\delta\lambda)({1-M\log\lmd})-\lambda}.
\endaligned
\de
\end{thm}

\proof
By Lemma \ref{baxmore}(ii) and the same approach used in Proposition \ref{atomerg},
for all $1<\lmd<e^{1/M}$ satisfying \rf{lmd1},
\be\lb{bnm}
\sup_{x_i\in A_1}\l|\l|\check{P}^n(x_i, \cdot)-\pi^*\r|\r|_\var
\leq \l(\pi^*(A_1)^{-1}-1\r)^{1/2}\lmd^{-n}=C_2\lmd^{-n}.
\de
Applying the first passage formula \rf{passage} to
the split chain $\check{\Phi}$, we obtain from
\rf{bnm} and Lemma \ref{baxmore}(ii) that for all $x_i\in\check{X}$,
\be\lb{361}
\aligned
||\check{P}^n(x_i, \cdot)-\pi^*||_\var
&\leq 2\check{\P}_{x_i}\{\check{\tau}_{A_1}\geq n+1\}
 +\sum_{m=1}^n\sup_{y_i\in A_1}||\check{P}^{n-m}(y_i, \cdot)
 -\pi^*||_\var\check{F}^m(x_i, A_1)\\
&\leq 2\check{\P}_{x_i}\{\check{\tau}_{A_1}\geq n+1\}
 +C_2\sum_{m=1}^n \lmd^{-(n-m)}\check{F}^m(x_i, A_1)\\
&\leq C_2 \lmd^{-n}+(2-C_2)^+\check{\P}_{x_i}
 \{\check{\tau}_{A_1}\geq n+1\}\\
&\q+C_2(1-\lmd^{-1})\lmd^{-n}\sum_{m=1}^n
 \lmd^{m}\check{\P}_{x_i}\{\check{\tau}_{A_1}\geq m\}\\
&\leq C_2 \lmd^{-n}
 +(2-C_2)^+\lmd^{-(n+1)}\sup_{x_i\in\check{X}}
 \check{\E}_{x_i}[\lmd^{\check{\tau}_{A_1}}]\\
&\q+C_2(1-\lmd^{-1})\lmd^{-n}\sum_{m=1}^n
 \sup_{x_i\in\check{X}}\check{\E}_{x_i}[\lmd^{\check{\tau}_{A_1}}]\\
&\leq C_2 \lmd^{-n}+(2-C_2)^+\lmd^{-(n+1)}M_2
 +C_2(1-\lmd^{-1})\lmd^{-n}M_2 n.\nnb\\
\endaligned
\de
By this inequality and noting that
$||P^n(x, \cdot)-\pi||_\var\leq||\check{P}^n(x_i, \cdot)-\pi^*||_\var$,
the desired assertion holds.
\qed

\bg{thm}
Under assumptions of Theorem \ref{minerg},
we have for all $1<\lmd<e^{1/M}$ satisfying \rf{lmd1},
$$
||\tld P^n-P^n||
\leq\fr{\lmd}{\lmd-1}\l[\l(D_2+\fr{E_2}{\lmd-1}\r)
-\l(D_2+\fr{E_2}{\lmd-1}+E_2n\r)\lmd^{-n}\r]||\tld P-P||,
$$
and
$$
||\tld\pi-\pi||_\var\leq\fr{\lmd}{\lmd-1}\l(D_2+\fr{E_2}{\lmd-1}\r)
||\tld P-P||,
$$
where $D_2$ and $E_2$ are defined in \rf{define-b}.
\end{thm}

%%%%%%%%%%%%%%%%%%%%%%%%%%%%%%%%%%%%%%%%%%%%%%%%%%%%%%%%%%%%%%%%%%%%%%%%%%%%%%%

\section{Uniform ergodicity and perturbation bounds for reversible Markov chains}\label{fr}

Unlike continuous-time Markov processes,
discrete-time Markov chains may not
be non-negative definite, this can cause troubles in our study.
For the description of this problem, see e.g. \ct{Bax05, MYH10, MS13, ST89}.
So more efforts need to be made to deal with general reversible Markov chains.
The method used here is to first investigate the two-skeleton chain with
the transition kernel $P^2$, and then transfer to $P$.

Let $\bar\Phi=\{\bar\Phi_n: n\in\Z\}$ be the Markov chain
with transition kernel $\bar P=P^2$.
It is obvious that $\bar\Phi$ is also reversible with respect to $\pi$
and always non-negative definite.
If $A$ is an atom for $P$, then for $x\in A$ and $B\in\mathscr{B}(X)$,
$$\bar P(x, B)=\nu(A)\nu(B)+\int_{A^c}P(y, B)\nu(d y),$$
i.e. $A$ is also an atom for $\bar P$
with the probability measure
$\bar\nu(\cdot):=\nu(A)\nu(\cdot)+\int_{A^c}P(y, \cdot)\nu(d y)$.
If the set $A$ satisfies
the minorization condition \rf{mincondition}, then
\be\lb{min2}
\bar P(x, B)\geq\int_A P(x, dy)P(y, B)
\geq\dlt^2\nu(A)\nu(B),\q x\in A, B\in\mathscr{B}(X).
\de
That is,
$A$ also satisfies the minorization condition
for $\bar{P}$ with the constant $\bar\dlt:=\dlt^2\nu(A)$
and the same probability measure $\nu$.

Let $\bar\tau_A=\inf\{n\geq1: \bar\Phi\in A\}$
be the first return time to $A$ for $\bar\Phi$,
and denoted by
$\bar F^n(x, A)=\P_x\{\bar\tau_A=n\}$
the distribution of $\bar\tau_A$.
For $s\geq0$, let
$$\bar F_{x A}(s)=\sum_{n=1}^{\ift}s^{2n}\bar F^n(x, A),$$
$$F^{(0)}_{x A}(s)=\sum_{n=1}^{\ift}s^{2n}F^{2n}(x, A),$$
and
$$F^{(1)}_{x A}(s)=\sum_{n=1}^{\ift}s^{2n-1}F^{2n-1}(x, A).$$
The next result shows the relationship
for the geometric moments of the first return times
$\tau_A$ and $\bar\tau_A$, which will be crucial for our method.
The related result for countable Markov chains
can be found in \ct[Proposition 2.1]{MYH10}.

\bg{lem}\lb{moments}
Let $A\in\mathscr{B}(X)$. Assume that
$
\sup_{x\in A}\sum_{n=1}^{\ift}F^{2n}(x, A)<1.
$
Then for $0\leq s\leq1$,
\be\lb{moment-relation}
\bar F_{x A}(s)\leq F^{(0)}_{x A}(s)+F^{(1)}_{x A}(s)\cdot
\sup_{y\in A}F^{(1)}_{y A}(s)\cdot\l[1-\sup_{y\in A}F^{(0)}_{y A}(s)\r]^{-1}.
\de
\end{lem}
\bg{proof}
Let the events $A_{n, 1}, \cdots, A_{n, n}$ be
$$A_{n, \ell}=\l\{\Phi_{m_1}\in A, \cdots, \Phi_{m_{\ell}}\in A~
\mbox{for}~\{m_1, \cdots, m_{\ell}\}\subset\{1, 3, \cdots, 2n-1\}~
\mbox{odd times}\r\}.$$
By the Markov property, we have
\be\lb{relation1}\aligned
\bar F^n(x, A)&=\P_x\l\{\bar\Phi_m\in A^c, 1\leq m<n, \bar\Phi_n\in A\r\}\\
&=\P_x\l\{\Phi_{2m}\in A^c, 1\leq m<n, \Phi_{2n}\in A\r\}\\
&=\P_x\l\{\Phi_{m}\in A^c, 1\leq m<2n, \Phi_{2n}\in A\r\}\\
&\q+\sum_{\ell=1}^{n}\P_x\l\{A_{n, \ell}, \Phi_{2m}\in A^c, 1\leq m<n,
  \Phi_{2n}\in A\r\}\\
&=F^{2 n}(x, A)+\sum_{m=1}^n\int_A {}_AP^{2 m-1}(x, d y){}_AP^{2 n-2 m+1}(y, A)+\cdots\\
&\q+\sum_{\scriptstyle{m_1+\cdots+m_{\ell}\leq n}\atop\scriptstyle{m_1, \cdots, m_{\ell}\geq 1}}\int_A\cdots\int_A {}_AP^{2m_1-1}(x, d y_1)
{}_AP^{2m_2}(y_1, d y_2)\cdots\\
&\q\q\q\q\q\q\q\q\q\q\q
\cdot{}_AP^{2 m_{\ell}}(y_{\ell-1}, d y_{\ell}){}_AP^{2n-2(m_1+\cdots+m_{\ell})+1}(y_{\ell}, A)\\
&\q+\cdots+\int_A\cdots\int_A {}_AP(x, d y_1)
{}_AP^{2}(y_1, d y_2)\cdots {}_AP^{2}(y_{n-1}, d y_n){}_AP(y_n, A).\\
\endaligned
\de
Noting that the summands in the right-hand side of \rf{relation1}
are multiple convolution, it follows that
\be\lb{2}\aligned
&\q\sum_{n=1}^{\ift}s^{2 n}\sum_{m=1}^n\int_A
{}_AP^{2 m-1}(x, d y){}_AP^{2 n-2 m+1}(y, A)\\
&\leq\sum_{m=1}^{\ift}s^{2 m-1}{}F^{2 m-1}(x, A)\cdot\sup_{y\in A}\sum_{n=1}^{\ift}s^{2 n-1}{}F^{2 n-1}(y, A)\\
&=F_{x A}^{(1)}(s)\cdot\sup_{y\in A}F_{y A}^{(1)}(s),\nnb\\
\endaligned
\de
where we use the fact ${}_AP^{2m-1}(x, A)=F^{2m-1}(x, A)$.
Similarly,
\be\lb{3}\aligned
&\q\sum_{n=1}^{\ift}s^{2 n}\sum_{\scriptstyle{m_1+\cdots+m_{\ell}\leq n}\atop\scriptstyle{m_1, \cdots, m_{\ell}\geq 1}}
\int_A\cdots\int_A {}_AP^{2m_1-1}(x, d y_1)
{}_AP^{2m_2}(y_1, d y_2)\cdots\\
&\q\q\q\q\q\q\q\q\q\q\q\q\q
\cdot{}_AP^{2 m_{\ell}}(y_{\ell-1}, d y_{\ell}){}_AP^{2n-2(m_1+\cdots+m_{\ell})+1}(y_{\ell}, A)\\
&\leq F_{x A}^{(1)}(s)\cdot\l[\sup_{y\in A}F_{y A}^{(0)}(s)\r]^{\ell-1}
\cdot\sup_{y\in A}F_{y A}^{(1)}(s).\nnb\\
\endaligned
\de
Multiplying both sides of \rf{relation1} by $s^{2 n}$
and making summation in $n$, we get
\be\lb{4}\aligned
\bar F_{x A}(s)&\leq\sum_{n=1}^{\ift}s^{2 n}F^{2n}(x, A)
+\sum_{\ell=1}^{\ift}F_{x A}^{(1)}(s)\cdot\l[\sup_{y\in A}F_{y A}^{(0)}(s)\r]^{\ell-1}\cdot\sup_{y\in A}F_{y A}^{(1)}(s)\\
&=F_{x A}^{(0)}(s)+F_{x A}^{(1)}(s)
\cdot\sup_{y\in A}F_{y A}^{(1)}(s)
\cdot\l[1-\sup_{y\in A}F_{y A}^{(0)}(s)\r]^{-1},\nnb\\
\endaligned
\de
where for $s\leq1$,
$$\sup_{y\in A}F_{y A}^{(0)}(s)\leq
\sup_{y\in A}\sum_{n=1}^{\ift}F^{2n}(y, A)<1.$$
This finishes the proof.
\end{proof}

Based on Lemma \ref{moments},
we derive the following Proposition \ref{buchong1}
and Theorem \ref{nondefinite-atom} when the state space contains
an accessible atom,
where Proposition \ref{buchong1} is the generalization of
\ct[Theorem 1.2]{MYH10}.

\bg{prop}\lb{buchong1}
For a reversible Markov chain,
assume that there exist some accessible atom
$A$ and some constant $\kp>1$ such that \rf{geo} holds.
Let
$$\rho=\sup\l\{s\leq\kp:\,
\sup_{x\in A}\sum_{n=1}^{\ift}s^{2n}F^{2n}(x, A)<1\r\}.$$
Then $r_0(P)\leq\rho^{-1}$. Moreover,
\be\lb{gen1}
\sup_{x\in A}||P^n(x, \cdot)-\pi||_\var\leq\l(\pi(A)^{-1}-1\r)^{1/2}\rho^{-n},\nnb
\de
and there exists a constant $C(x)<\ift$ such that
\be\lb{gen2}
||P^n(x, \cdot)-\pi||_\var\leq C(x)\rho^{-n}, \q\pi\mbox{-a.s.}~x\in X.\nnb
\de
\end{prop}

\bg{proof}
For the accessible atom $A$, it is obvious that
$$\sup_{x\in A}\sum_{n=1}^{\ift}F^{2n}(x, A)
=\sum_{n=1}^{\ift}F^{2n}(x, A)<1,\q x\in A.$$
Notice that although \rf{moment-relation} is proved
only for $0\leq s\leq1$, \rf{moment-relation} is true
for any $s$ such that
$$\sup_{x\in A}\sum_{n=1}^{\ift}s^{2n}F^{2n}(x, A)<1.$$
It follows from Lemma \ref{moments} that for all $s<\rho$,
$\sup_{x\in A}\E_x\l[s^{2\bar\tau_A}\r]<\ift$.
Then a similar proof of Proposition \ref{atomerg} implies that
$r_0(\bar P)\leq s^{-2}$ for all $s<\rho$,
so that $r_0(\bar P)\leq\rho^{-2}$.
Thus, $r_0(P)\leq\rho^{-1}$ and the desired results hold.
\end{proof}

\bg{thm}\lb{nondefinite-atom}
For a reversible Markov chain, assume that there exists
some accessible atom $A$ such that \rf{hitting} holds.
Let
$$\varrho=\sup\l\{s<e^{1/M}:\,
\sup_{x\in A}\sum_{n=1}^{\ift}s^{2n}F^{2n}(x, A)<1\r\}.$$
Then for all $1<\lmd<e^{1/M}$,
\be\lb{bus1a}
\aligned
\l|\l|P^n-\pi\r|\r|\leq
\left\{\begin{array}{ll}
F_1\varrho^{-n}+G_1\lmd^{-n},& \lmd\neq\varrho;\\
\l(J_1+K_1 n\r)\varrho^{-n},& \lmd=\varrho,\\
\end{array} \right.\nnb\\
\endaligned
\de
where
\be\lb{define-f}
\aligned
&F_1=C_1\l(1-\fr{\varrho-1}{\varrho-\lmd}M_1\r),\q
G_1=M_1\l((2-C_1)^+\lmd^{-1}+\fr{\varrho-1}{\varrho-\lmd}C_1\r),\\
&J_1=C_1+(2-C_1)^+\varrho^{-1}M_1,\q K_1=C_1(1-\varrho^{-1})M_1,\\
\endaligned
\de
and $C_1$ and $M_1$ are defined in \rf{changshu1}.
\end{thm}

\bg{proof}
According to Lemma \ref{hitmoment} and Lemma \ref{moments}, 
$\sup_{x\in X}\E_x[s^{2\bar\tau_A}]<\ift$ for all $s<\varrho$.
Hence $r_0(P)\leq\varrho^{-1}$, and then
$$
\sup_{x\in A}||P^n(x, \cdot)-\pi||_\var\leq\l(\pi(A)^{-1}-1\r)^{1/2}\varrho^{-n}.
$$
Combining this inequality with Lemmas \ref{renewal} and \ref{hitmoment},
for all $x\in X$ and $1<\lmd<e^{1/M}$,
\be\lb{sdfg}
\aligned
||P^n(x, \cdot)-\pi||_\var
&\leq 2\P_x\{\tau_A\geq n+1\}
+C_1\sum_{m=1}^n\varrho^{-(n-m)}\P_x\l\{\tau_A=m\r\}\\
&\leq C_1\varrho^{-n}+(2-C_1)^+\sup_{x\in X}\P_x\{\tau_A\geq n+1\}\\
&\q+C_1(1-\varrho^{-1})\varrho^{-n}
\sum_{m=1}^n\varrho^{m}\sup_{x\in X}\P_x\{\tau_A\geq m\}\\
&\leq C_1\varrho^{-n}+(2-C_1)^+\lmd^{-(n+1)}\sup_{x\in X}\E_x[\lmd^{\tau_A}]\\
&\q+C_1(1-\varrho^{-1})\varrho^{-n}
\sum_{m=1}^n\l(\varrho/\lmd\r)^{m}\sup_{x\in X}\E_x[\lmd^{\tau_A}]\\
&\leq  C_1\varrho^{-n}+(2-C_1)^+\lmd^{-(n+1)}M_1
 +C_1(1-\varrho^{-1})\varrho^{-n}M_1\sum_{m=1}^n\l(\varrho/\lmd\r)^m.\nnb\\
\endaligned
\de
Then the desired result holds by the above inequality
and simple calculations.
\end{proof}

For non-atomic case, through the same splitting techniques
of Section \ref{nonatom},
we can split the Markov chain $\bar\Phi$
to produce a new chain $\hat\Phi$
on $(\hat{X}, \mathscr{B}(\hat{X}))$
with transition kernel $\hat{P}$ and invariant probability
measure $\pi^{**}$, which contains an accessible atom $A_1$.
It is obvious from \rf{measure} and \rf{min2} that
$\pi^{**}(A_1)=\dlt^2\nu(A)\pi(A)$.
Let $\hat{\E}_{x_i}$ be the expectation for $\hat\Phi$
started with $\hat{\Phi}_0=x_i$, and denote by
$\hat{\tau}_{A_1}=\inf\{n\geq1: \hat{\Phi}_n\in A_1\}$
the first return time to $A_1$ for $\hat\Phi$.
%and
%$$\hat{F}^n(x_i, A_1)=\hat{\P}_{x_i}\{\hat{\tau}_{A_1}=n\},
%\q n\in\N\cup\{\ift\}$$
%the distribution of $\hat{\tau}_{A_1}$.

\bg{prop}
For a reversible Markov chain, assume that there exist some set
$A\in\mathscr{B}^+(X)$ and some constant $\kp>1$ such that
\rf{mincondition} and \rf{geo} hold, and
\be\lb{cc}
\sup_{x\in A}\sum_{n=1}^{\ift}\kp^{2n}F^{2n}(x, A)\leq\tht
\de
for some $\tht<1$.
Then there exists a constant $C(x)<\ift$ such that
$$
||P^n(x, \cdot)-\pi||_\var\leq C(x)\Gm^{-n}, \q\pi\mbox{-a.s.}~x\in X,
$$
where $\Gm=\kp\wedge(1-\dlt^2\nu(A))^{-1/\gm}$ and
$$\gm=\l(\log\fr{L^2/(1-\tht)-\dlt^2\nu(A)\kp^2}
{1-\dlt^2\nu(A)}\r)\big/\l(\log\kp\r).$$
\end{prop}

\bg{proof}
Under \rf{geo} and \rf{cc},
we get by Lemma \ref{moments} that
\be\lb{zuj}
\sup_{x\in A}\E_x\l[\kp^{2\bar{\tau}_A}\r]\leq \fr{L^2}{1-\tht}.\nnb
\de
Applying the techniques used in Lemma \ref{baxmore}(i)
to the chain $\bar\Phi$, we obtain from \rf{min2} and the above inequality that
for all $1<\lmd<\kp\wedge(1-\dlt^2\nu(A))^{-1/\gm}$,
$$\hat{\E}_{A_1}[\lmd^{2\hat\tau_{A_1}}]<\ift.$$
Thus, a similar proof as that of Proposition \ref{minerg-prop} yields that
$$
||P^{2n}(x, \cdot)-\pi||_\var\leq C(x)\Gm^{-2n}, \q\pi\mbox{-a.s.}~x\in X.
$$
From this and noting that
\be\lb{yasuo}
||P^{2n+1}(x, \cdot)-\pi||_{\var}\leq||P^{2n}(x, \cdot)-\pi||_{\var},
\de
the desired assertion holds.
\end{proof}

\bg{thm}\lb{nondefinite}
For a reversible Markov chain, assume that there exists a set
$A\in\mathscr{B}^+(X)$ such that \rf{mincondition} and \rf{hitting} hold.
Then for all $1<\lmd<e^{1/M}$ satisfying
\be\lb{tht}
\sup_{x\in A}\sum_{n=1}^{\ift}\lmd^{2n}F^{2n}(x, A)\leq\vartheta,
\de
and
\be\lb{thtt}
\lmd^2<(1-\vartheta)(1-M\log\lmd)^2(1+\dlt^2\nu(A)\lmd^2)
\de
for some $\vartheta<1$,
\be\lb{ct6}
\l|\l|P^n-\pi\r|\r|\leq\l(D_3+E_3 n\r)\lmd^{-n},\nnb
\de
where
\be\lb{define-3}
\aligned
&D_3=C_3\lmd+(2-C_3)^+\lmd^{-1}M_3,\q
E_3=C_3(\lmd-\lmd^{-1})M_3/2,\\
&C_3=\l((\dlt^2\nu(A)\pi(A))^{-1}-1\r)^{1/2},\;
M_3=\fr{\dlt^2\nu(A)\lmd^2}{(1-\vartheta)
(1-M\log\lmd)^2(1+\dlt^2\nu(A)\lmd^2)-\lmd^2}.\\
\endaligned
\de
\end{thm}

\bg{proof}
According to Lemma \ref{hitmoment} and Lemma \ref{moments},
for all $1<\lmd<e^{1/M}$ satisfying \rf{tht},
\be\lb{zui1}
\sup_{x\in X}\E_x\l[\lmd^{2\bar{\tau}_A}\r]\leq \fr{\lmd^2}{(1-\vartheta)(1-M\log\lmd)^2}=:\bar M.
\de
Applying the techniques used in Lemma \ref{baxmore}(ii)
to the chain $\bar\Phi$, we get from \rf{min2} and \rf{zui1} that
for all $1<\lmd<e^{1/M}$ satisfying \rf{tht} and \rf{thtt},
$$\sup_{x_i\in\hat{X}}\hat{\E}_{x_i}[\lmd^{2\hat{\tau}_{A_1}}]
\leq\fr{\dlt^2\nu(A)\bar M}
{1-\l(\bar M-\dlt^2\nu(A)\lmd^2\r)}=M_3.$$
Then by a similar proof as that of Theorem \ref{minerg},
\be\lb{contarction1}
||P^{2n}(x, \cdot)-\pi||_{\var}
\leq C_3 \lmd^{-2n}+(2-C_3)^+\lmd^{-2(n+1)}M_3
+C_3(1-\lmd^{-2})\lmd^{-2n}M_3 n.\nnb
\de
Thus, the desired result holds by the above inequality and \rf{yasuo}.
\end{proof}

Combining Theorems \ref{sensitivity} with
\ref{nondefinite-atom} or \ref{nondefinite},
we obtain immediately the following perturbation bounds.

\bg{thm}
$(i)$
    Under assumptions of Theorem \ref{nondefinite-atom},
    we have for all $1<\lmd<e^{1/M}$,
\be\lb{zuihou1}
    \aligned
    ||\tld P^n-P^n||\leq
    \left\{\begin{array}{ll}
    \l[\fr{\varrho}{\varrho-1}F_1(1-\varrho^{-n})+
    \fr{\lmd}{\lmd-1}G_1(1-\lmd^{-n})\r]
    ||\tld P-P||,& \lmd\not=\varrho;\\
    \fr{\varrho}{\varrho-1}\l[\l(J_1+\fr{K_1}{\varrho-1}\r)
    -\l(J_1+\fr{K_1}{\varrho-1}+K_1n\r)\varrho^{-n}\r]||\tld P-P||,& \lmd=\varrho,\\
    \end{array} \right.\nnb\\
    \endaligned
\de
and
\be\lb{zuihou2}
    \aligned
    ||\tld\pi-\pi||_\var\leq
    \left\{\begin{array}{ll}
    \l[\fr{\varrho}{\varrho-1}F_1+\fr{\lmd}{\lmd-1}G_1\r]||\tld P-P||,& \lmd\not=\varrho;\\
    \fr{\varrho}{\varrho-1}\l(J_1+\fr{K_1}{\varrho-1}\r)||\tld P-P||,& \lmd=\varrho,\\
    \end{array} \right.\nnb\\
    \endaligned
\de
where $F_1$, $G_1$, $J_1$ and $K_1$ are defined in \rf{define-f}.

$(ii)$
Under assumptions of Theorem \ref{nondefinite},
we have for all $1<\lmd<e^{1/M}$ satisfying \rf{tht} and \rf{thtt},
\be\lb{zuihou3}
\aligned
||\tld P^n-P^n||
\leq\fr{\lmd}{\lmd-1}\l[\l(D_3+\fr{E_3}{\lmd-1}\r)
-\l(D_3+\fr{E_3}{\lmd-1}+E_3 n\r)\lmd^{-n}\r]||\tld P-P||,\nnb\\
\endaligned
\de
and
\be\lb{zuihou4}
\aligned
||\tld\pi-\pi||_\var\leq\fr{\lmd}{\lmd-1}\l(D_3+\fr{E_3}{\lmd-1}\r)
||\tld P-P||,\nnb\\
\endaligned
\de
where $D_3$ and $E_3$ are defined in \rf{define-3}.
\end{thm}

%%%%%%%%%%%%%%%%%%%%%%%%%%%%%%%%%%%%%%%%%%%%%%%%%%%%%%%%%%%%%%%%%%%%%%%%%%%%%%%
\section{Perturbation bounds for general Markov chains}\label{gen}

In the section, by using a result in \ct{ak},
we present a different bound for $||\tld\pi-\pi||_{\var}$
via the uniform moments of the first hitting times.
The advantage of the estimate is that it works for
general (non-reversible) Markov chains.

Let $V\geq1$ be a measurable function.
Recall that the $V$-norm distance between two transition kernels
$\tld P$ and $P$ is defined as
$$|||\tld P-P|||_{V}=\sup_{x\in X}
\fr{||\tld P(x, \cdot)-P(x, \cdot)||_{V}}{V(x)},$$
where for any signed measure $\mu$ on $\mathscr{B}(X)$,
$||\mu||_{V}=\sup_{|f|\leq V}|\mu(f)|$.
For general Markov chains, the following two conditions
are used in the sequel:

\noindent
(A1)\; There exist a probability measure $\nu$
      on $\mathscr{B}(X)$ and a bounded non-negative function

      \q\;$h$ with $\pi(h)>0$ and $\nu(h)>0$ such that
      $$T(x, B):=P(x, B)-h(x)\nu(B)\geq0,\q x\in X, B\in\mathscr{B}(X).$$

\noindent
(A2)\; There exist some function $V\geq1$
      and some constant $\rho<1$ such that
      $$|||P|||_{V}<\ift,\q
      T V(x)\leq\rho V(x),\q x\in X.$$
Clearly, (A1) can be seen as a more general minorization condition,
for more details the interested readers should consult \ct{kar1, Num84}.
Under (A1) and (A2), we have the next lemma by \ct{ak}
or \ct[Theorem 3.8 and Remark 2.5]{kar1}.

\bg{lem}\lb{noncond}
For a general Markov chain, assume that $(A1)$ and $(A2)$ hold.
Then for $|||\tld P-P|||_{V}<(1-\rho)/c$, we have
\be\lb{baibai}
||\tld\pi-\pi||_{V}\leq c||\pi||_{V}
\l(1-\rho-c\,|||\tld P-P|||_{V}\r)^{-1}|||\tld P-P|||_{V},\nnb
\de
where
$c=1+||\pi||_{V}/\inf_{x\in X}V(x)$.
\end{lem}

Lemma \ref{noncond} enables one to derive
the following perturbation bound.

\bg{thm}\lb{general MC}
For a general Markov chain,
assume that there exists a set $A\in\mathscr{B}^+(X)$
such that \rf{hitting} holds.
Then for all $1<\lmd<e^{1/M}$,
$$||\tld\pi-\pi||_{\var}\leq M_0^2(1+M_0)
\l(1-\lmd^{-1}-M_0(1+M_0)||\tld P-P||\r)^{-1}||\tld P-P||$$
provided
$||\tld P-P||<(1-\lmd^{-1})/(M_0+M_0^2)$, where
$M_0=(1-M\log\lmd)^{-1}$.
\end{thm}

\bg{proof}
First, we check conditions (A1) and (A2) hold.
Let $h(x)=1_A(x)$ and $\nu=P(x, \cdot)$ for $x\in X$.
Then (A1) holds obviously.
Set $V(x)=\E_{x}\l[\lmd^{\sgm_A}\r]$ for $x\in X$.
Then we have $1\leq V\leq M_0$ by \rf{m0}, and
$$|||P|||_{V}=\sup_{x\in X}\fr{P V(x)}{V(x)}\leq M_0.$$
From \ct[Corollary 2.8]{MS14},
$$T V(x)=P V(x)=\lmd^{-1}V(x),\q x\in A^c,$$
and
$$T V(x)=0\leq\lmd^{-1}V(x),\q x\in A.$$
That is, condition (A2) is fulfilled with
$V(x)=\E_{x}\l[\lmd^{\sgm_A}\r]$ and $\rho=\lmd^{-1}$.

Next, to apply Lemma \ref{noncond},
we shall compute $||\pi||_{V}$ and $c$.
Since $1\leq V\leq M_0$, we have
$||\pi||_{V}\leq M_0$ and $c\leq 1+M_0$.
Hence for $||\tld P-P||<(1-\lmd^{-1})/(M_0+M_0^2)$,
$$|||\tld P-P|||_{V}\leq M_0 ||\tld P-P||
<\fr{1-\lmd^{-1}}{1+M_0}\leq\fr{1-\rho}{c},
$$
so we get
\be\lb{666}
\aligned
&\q||\tld\pi-\pi||_\var
\leq||\tld\pi-\pi||_V\\
&\leq c||\pi||_{V}
\l(1-\rho-c\,|||\tld P-P|||_{V}\r)^{-1}|||\tld P-P|||_{V}\\
&\leq M_0^2(1+M_0)
\l(1-\lmd^{-1}-M_0(1+M_0)||\tld P-P||\r)^{-1}||\tld P-P||,\nnb\\
\endaligned
\de
which is the desired assertion.
\end{proof}

%%%%%%%%%%%%%%%%%%%%%%%%%%%%%%%%%%%%%%%%%%%%%%%%%%%%%%%%%%%%%%%%%%%%%%%%%%%%%%%
\bigskip
\bigskip
\noindent\textbf{Acknowledgements}
This work is supported by the National Natural Science
Foundation of China (Grant Nos. 11771047, 11501576).

\end{document}